\documentclass{amsart}
\usepackage{amsfonts}

\newtheorem{theorem}{Theorem}[section]
\newtheorem{lemma}[theorem]{Lemma}
\newtheorem{proposition}[theorem]{Proposition}
\newtheorem{corollary}[theorem]{Corollary}
\numberwithin{equation}{section}

\newcommand{\0}[2]{\langle #1,#2\rangle}
\newcommand{\1}[2]{\langle #1,#2\rangle_1}
\newcommand{\2}[2]{\langle #1,#2\rangle_2}
\newcommand{\re}{\operatorname{Re}}

\author{Minghua Lin}
\address{
Department of Combinatorics and Optimization\\ University of Waterloo\\
Waterloo, Canada}
\email{mlin87@ymail.com}

\author{Gord Sinnamon}
\address{
Department of Mathematics\\
University of Western Ontario\\
London, Canada}
\email{sinnamon@uwo.ca}
\thanks{Supported by the Natural Sciences and Engineering Research Council of Canada}


\begin{document}

\title[The Generalized Wielandt Inequality]{The Generalized Wielandt Inequality in Inner Product Spaces}

\begin{abstract}
A new inequality between angles in inner product spaces is formulated and proved. It leads directly to a concise statement and proof of the generalized Wielandt inequality, including a simple description of all cases of equality. As a consequence, several recent results in matrix analysis and inner product spaces are improved.
\end{abstract}
\maketitle

\section{Introduction}\label{intro}
The Wielandt and generalized Wielandt inequalities control how much angles can change under a given invertible matrix transformation of $\mathbb C^n$. The control is given in terms of the condition number of the matrix. Wielandt, in \cite{W}, gave a bound on the resulting angles when orthogonal complex lines are transformed. Subsequently, Bauer and Householder, in \cite{BH}, extended the inequality to include arbitrary starting angles. These basic inequalities of matrix analysis were introduced to give bounds on convergence rates of iterative projection methods but have found a variety of applications in numerical methods, especially eigenvalue estimation. They are also applied in multivariate analysis, where angles between vectors correspond to statistical correlation. See, for example, \cite{BH}, \cite{E}, \cite{ET}, \cite{HJ} and  \cite{H}. There are also matrix-valued versions of the inequality that are receiving attention, especially in the context of statistical analysis. See \cite{BD}, \cite{LLP}, \cite{WI}, and \cite{ZZ}. 

The condition number of an invertible matrix $A$ is $\kappa(A)=\|A\|\|A^{-1}\|$, where $\|\cdot\|$ denotes the operator norm. If $A$ is positive definite and Hermitian, $\kappa(A)$ is easily seen to be the ratio of the largest and smallest eigenvalues of $A$. The following statement of the generalized Wielandt inequality is taken from \cite{H}.

\begin{theorem}\label{gW} Let $A$ be an invertible $n\times n$ matrix. If $x,y\in \mathbb C^n$ and $\Phi,\Psi\in[0,\pi/2]$ satisfy 
\[
|y^*x|\le\|x\|\|y\|\cos\Phi\quad\mbox{and}\quad\cot(\Psi/2)=\kappa(A)\cot(\Phi/2),
\] 
then 
\[
|(Ay)^*(Ax)|\le\|Ax\|\|Ay\|\cos\Psi.
\]
\end{theorem}

The generalized Wielandt inequality can be difficult to apply for several reasons. First, despite having various equivalent formulations, the inequality seems always to be expressed in ways that hide the natural symmetry coming from the invertible transformation involved. Next, the conditions for equality are known, see \cite{K}, but are unwieldy and hard to apply. Finally, the angles involved are angles between complex lines rather than between individual vectors.

Although the last point seems minor, we found it to be the key to a symmetric formulation and a simple description of the cases of equality. In Theorem \ref{main} and its matrix analytic counterpart, Theorem \ref{Mmain}, we present a new inequality that gives sharp upper and lower bounds for the angle between a pair of transformed vectors. The conditions for equality are simple and easy to apply. This new inequality relates angles between vectors rather than between complex lines but it immediately implies a result for angles between complex lines that is equivalent to the generalized Wielandt inequality. Moreover, this version of the generalized Wielandt inequality retains the simple form of the new inequality and (most of) the simplicity of its conditions for equality.

In Section \ref{mainResults} we work in the context of an arbitrary real or complex vector space having two inner products. This approach preserves symmetry by avoiding the distinction between angles before and after a fixed transformation. Also, the main result is not restricted to $\mathbb C^n$ but holds for vectors in infinite-dimensional spaces. As an application of the unrestricted result, we improve a metric space inequality from \cite{D}. The main results are then formulated in the language of matrix analysis in Section \ref{matrices}, and we apply them to improve inequalities from \cite{Yeh} and \cite{LS}, and to settle a conjecture from \cite{Yan}.

To begin, a short discussion of angles in inner product spaces is in order. In a real inner product space $(V,\0\cdot\cdot)$ the angle $\theta=\theta(u,v)$ between two non-zero vectors is defined by, $0\le\theta\le\pi$ and 
\[
\cos\theta=\frac{\0uv}{\|u\|\|v\|}.
\]
Here $\|u\|=\sqrt{\0uu}$ is the norm induced by the inner product. The angle between subsets $S$ and $T$ of $V$ is the infimum of the angles between non-zero elements of $S$ and $T$, so 
\[
\Theta(S,T)=\inf\{\theta(u,v):0\ne u\in S, 0\ne v\in T\}.
\]
With this definition it is easy to check that the angle $\Theta=\Theta(\mathbb R u,\mathbb R v)$ between the lines $\mathbb R u$ and $\mathbb R v$ satisfies $0\le\Theta\le\pi/2$ and
\[
\cos\Theta=\frac{|\0uv|}{\|u\|\|v\|}.
\]

A complex inner product space $(V,\0\cdot\cdot)$ may be viewed as the real inner product space $(V_\mathbb R,\re\0\cdot\cdot)$ where $V_\mathbb R=V$ with the scalars restricted to $\mathbb R$. Since $\re\0vv=\0vv$ for all $v\in V$, lengths in $V$ are preserved and therefore so are angles. Thus, this real inner product is used to define the angle $\theta$ between the vectors $u$ and $v$, and a computation gives the formula for the angle $\Theta$ between the complex lines $\mathbb C u$ and $\mathbb C v$. We have, 
\[
\cos\theta=\frac{\re\0uv}{\|u\|\|v\|}\quad\mbox{and}\quad\cos\Theta=\frac{|\0uv|}{\|u\|\|v\|}.
\]
The second formula is often used as a definition of the angle between vectors $u$ and $v$ in a complex inner product space. (Angles defined this way do not determine angles in triangles correctly but they have the advantage that complex orthogonality, namely $\0uv=0$, is equivalent to the angle between $u$ and $v$ being $\pi/2$.)

We will make use of the simple observation that if $|\alpha|=1$, then
\begin{equation}\label{alpha}
\Theta(\mathbb C u,\mathbb C v)=\theta(\alpha u,v)\quad\mbox{if and only if}\quad|\0uv|=\alpha\0uv.
\end{equation}
(Note that our inner products are taken to be linear in the first variable.) The above observation remains valid for $\Theta(\mathbb R u,\mathbb R v)$ in a real inner product space, where $\alpha=\pm1$.

\section{Main results}\label{mainResults}
Suppose $V$ is a non-trivial real or complex vector space. Let $\1\cdot\cdot$ and $\2\cdot\cdot$ be inner products on $V$ and define  $m$, $V_m$, $M$,  $V_M$, $E_1$ and $E_2$ by,
\begin{equation}\label{defs}
\left\{\begin{gathered}
m=\inf_{0\ne v\in V} \|v\|_2\left/\|v\|_1\right.,\quad V_m=\{v\in V:\|v\|_2=m\|v\|_1\},\\
M=\sup_{0\ne v\in V} \|v\|_2\left/\|v\|_1\right.,\quad V_M=\{v\in V:\|v\|_2=M\|v\|_1\},\\
E=E_j=\left\{(u,v):\frac u{\|u\|_j}+\frac v{\|v\|_j}\in V_m,\frac u{\|u\|_j}-\frac v{\|v\|_j}\in V_M\right\},
\end{gathered}\right.
\end{equation}
for $j=1,2$. Here, as usual, $\|v\|_1=\sqrt{\1vv}$ and $\|v\|_2=\sqrt{\2vv}$. We anticipate the result of Corollary \ref{E} in the definition of $E$ above.

Evidently $0\le m\le M\le\infty$, $0\in V_m$ and $0\in V_M$. (The convention $0\cdot\infty=0$ ensures that $0\in V_M$ when $M=\infty$.)  A standard compactness argument shows that if $V$ is finite dimensional then $0<m\le M<\infty$ and $V_m\ne\{0\}\ne V_M$. If $m=M$ then $V_m=V_M=V$ and, by polarization, $\2uv=m^2\1uv$ for all $u,v\in V$. 

\begin{lemma}\label{orthogonal} Let $V$ be a real vector space equipped with inner products $\1\cdot\cdot$ and $\2\cdot\cdot$. Make the definitions (\ref{defs}). If $m<M$, then $V_m$ and $V_M$ are subspaces and the two are mutually orthogonal with respect to both inner products.
\end{lemma}
Proof. Suppose $u$ is a non-zero vector in $V_m$ and $v\in V$ is not a multiple of $u$. Then
\[
f(t)=\frac{\|u+tv\|_2^2}{\|u+tv\|_1^2}=\frac{\2uu+2t\2uv+t^2\2vv}{\1uu+2t\1uv+t^2\1vv}
\]
is defined and differentiable for $t\in \mathbb R$. Since $f$ achieves its minimum value at $t=0$, $f'(0)=0$. That is, $\2uv\1uu=\2uu\1uv$. Thus, for all $u\in V_m$ and all $v\in V$, 
\[
\2uv=m^2\1uv.
\]
(The excluded case, $u=0$ or $v$ a multiple of $u$, is easily verified.)
It follows that if $v\in V_m$ then $f$ is the constant function with value $m^2$. In particular, $f(1)=m^2$, so $u+v\in V_m$. Since it is clearly closed under scalar multiplication, $V_m$ is a subspace. 

Repeating the argument for $V_M$ shows that it, too, is a subspace and that for all $v\in V_M$ and $u\in V$,
\[
\2uv=M^2\1uv.
\]
If $u\in V_m$ and $v\in V_M$ then $m^2\1uv=\2uv=M^2\1uv$ and hence $\1uv=\2uv=0$. Thus $u$ and $v$ are orthogonal with respect to both inner products. This completes the proof.\qed\vskip 2ex

\begin{corollary}\label{2basis} Let $V$ be a real vector space equipped with inner products $\1\cdot\cdot$ and $\2\cdot\cdot$. Make the definitions (\ref{defs}). If $V$ is two-dimensional, then there is a basis of $V$ that is orthogonal with respect to both inner products.
\end{corollary}
Proof. If $m=M$ then the two inner products are multiples of each other and any orthogonal basis will do. Otherwise, let $0\ne b\in V_m$ and $0\ne B\in V_M$. Then $\{b,B\}$ is the desired basis.\qed\vskip 2ex

The next result justifies the use of $E$ to denote either $E_1$ or $E_2$.
\begin{corollary}\label{E} Let $V$ be a real vector space equipped with inner products $\1\cdot\cdot$ and $\2\cdot\cdot$. Make the definitions (\ref{defs}). Then $E_1=E_2$.
\end{corollary}
Proof. By symmetry it is enough to show that $E_1\subseteq E_2$. For $(u,v)\in E_1$, let
\[
w=\frac u{\|u\|_1}+\frac v{\|v\|_1}\in V_m\quad\mbox{and}\quad W=\frac u{\|u\|_1}-\frac v{\|v\|_1}\in V_M.
\]
By Lemma \ref{orthogonal}, $w$ and $W$ are orthogonal with respect to $\2\cdot\cdot$, so
\[
\|u\|_2^2\left/\|u\|_1^2\right.=\tfrac14\|w+W\|_2^2=\tfrac14(\|w\|_2^2+\|W\|_2^2)=\tfrac14\|w-W\|_2^2=\|v\|_2^2\left/\|v\|_1^2\right..
\]
Thus 
\[
\frac u{\|u\|_2}+\frac v{\|v\|_2}=\frac{\|u\|_1}{\|u\|_2}w\in V_m\quad\mbox{and}\quad \frac u{\|u\|_2}-\frac v{\|v\|_2}=\frac{\|u\|_1}{\|u\|_2}W\in V_m
\]
and so $(u,v)\in E_2$.\qed\vskip 2ex

Having two inner products, the space $V$ has two differing notions of the angle between vectors. Our main result provides a comparison between these angles in terms of the quantities $m$ and $M$ defined in (\ref{defs}).

\begin{theorem}\label{main} Let $V$ be a real or complex vector space equipped with inner products $\1\cdot\cdot$ and $\2\cdot\cdot$. Make the definitions (\ref{defs}). For independent vectors $u$ and $v$ in $V$ let $\varphi$ and $\psi$ be defined by, $0\le\varphi\le\pi$, $0\le\psi\le\pi$,
\[
\cos\varphi=\frac{\re\1uv}{\|u\|_1\|v\|_1}\quad\mbox{and}\quad\cos\psi=\frac{\re\2uv}{\|u\|_2\|v\|_2}.
\]
Then 
\begin{equation}\label{TAN}
(m/M)\tan(\varphi/2)\le\tan(\psi/2)\le (M/m)\tan(\varphi/2).
\end{equation}
Equality holds in the right-hand inequality if and only if $(u,v)\in E$. Equality holds in the left-hand inequality if and only if $(u,-v)\in E$.
\end{theorem}
Proof.  First consider the case that $V$ is a real vector space. Note that the assumption of independence ensures $0<\varphi<\pi$ and  $0<\psi<\pi$.

By Corollary \ref{2basis}, the span of $u$ and $v$ has a basis $\{b,B\}$ that is orthogonal with respect to both inner products. Without loss of generality we may assume that $\|b\|_1=\|B\|_1=1$. For notational convenience, set $n=\|b\|_2$ and $N=\|B\|_2$ and suppose, by interchanging $b$ and $B$ if necessary, that $n\le N$.  Note that the definitions of $m$ and $M$ ensure that $m\le n$ and $N\le M$. Write $u=u_bb+u_BB$ and $v=v_bb+v_BB$ for some real numbers $u_b$, $u_B$, $v_b$, and $v_B$. In terms of these coordinates we have,
\begin{align*}
\|u\|_1^2\|v\|_1^2\sin^2\varphi
&= \|u\|_1^2\|v\|_1^2-\1uv^2\\
&=(u_b^2+u_B^2)(v_b^2+v_B^2)-(u_bv_b+u_Bv_B)^2\\
&=(u_bv_B-u_Bv_b)^2
\end{align*}
and
\begin{align*}
\|u\|_2^2\|v\|_2^2\sin^2\psi
&= \|u\|_2^2\|v\|_2^2-\2uv^2\\
&=(n^2u_b^2+N^2u_B^2)(n^2v_b^2+N^2v_B^2)-(n^2u_bv_b+N^2u_Bv_B)^2\\
&=n^2N^2(u_bv_B-u_Bv_b)^2.
\end{align*}
Thus,
\begin{equation}\label{sin}
\|u\|_2\|v\|_2\sin\psi=nN\|u\|_1\|v\|_1\sin\varphi.
\end{equation}

The derivative of
\[
g(x) = (u_b^2+xu_B^2)^{1/2}(v_b^2+xv_B^2)^{1/2}+(u_bv_b+xu_Bv_B)
\]
is
\[
g'(x)=\frac12\left(u_B\left(\frac{v_b^2+xv_B^2}{u_b^2+xu_B^2}\right)^{1/4}+v_B\left(\frac{u_b^2+xu_B^2}{v_b^2+xv_B^2}\right)^{1/4}\right)^2\ge0,
\]
so $g(1)\le g(N^2/n^2)$.  Multiplying both sides of this by $n^2$ gives,
\begin{equation}\label{n2cos}
n^2\|u\|_1\|v\|_1(1+\cos\varphi)\le\|u\|_2\|v\|_2(1+\cos\psi).
\end{equation}
Combining (\ref{sin}) and (\ref{n2cos}) gives,
\begin{equation}\label{nN}
\tan(\psi/2)=\frac{\sin\psi}{(1+\cos\psi)}
\le\frac{nN\sin\varphi}{n^2(1+\cos\varphi)}=(N/n)\tan(\varphi/2),
\end{equation}
with equality if and only if $g'(x)=0$ for $x\in (1,N^2/n^2)$. Since $m\le n\le N\le M$, (\ref{nN}) proves the right-hand inequality of (\ref{TAN}). 

If equality holds in the right-hand inequality of (\ref{TAN}), then equality holds in (\ref{nN}) and $n=m$, $N=M$, $b\in V_m$, and $B\in V_M$. If $m=M$ then $V_m=V_M=V$ and $\varphi=\psi$ so the last two statements of the theorem are trivial. Otherwise, equality in (\ref{nN}) implies that $g'$ is zero on the non-trivial interval $(1,M^2/m^2)$. That is, 
\[
u_B\left(\frac{v_b^2+xv_B^2}{u_b^2+xu_B^2}\right)^{1/4}+v_B\left(\frac{u_b^2+xu_B^2}{v_b^2+xv_B^2}\right)^{1/4}=0
\]
and hence $u_B^2v_b^2=v_B^2u_b^2$. Since $u$ and $v$ are independent, both $u_B$ and $v_B$ are non-zero, they have opposite signs, and $u_Bv_b=- v_Bu_b$.  Therefore,
\begin{align*}
&\frac u{\|u\|_1}+\frac v{\|v\|_1}
=\frac{u_bb+u_BB}{\sqrt{u_b^2+u_B^2}}+\frac{v_bb+v_BB}{\sqrt{v_b^2+v_B^2}}\\
&=\pm\left(\frac{(u_b/u_B)b+B}{\sqrt{(u_b/u_B)^2+1}}-\frac{(v_b/v_B)b+B}{\sqrt{(v_b/v_B)^2+1}}\right)
=\pm\frac{2(u_b/u_B)b}{\sqrt{(u_b/u_B)^2+1}}\in V_m
\end{align*}
and
\begin{align*}
&\frac u{\|u\|_1}-\frac v{\|v\|_1}
=\frac{u_bb+u_BB}{\sqrt{u_b^2+u_B^2}}-\frac{v_bb+v_BB}{\sqrt{v_b^2+v_B^2}}\\
&=\pm\left(\frac{(u_b/u_B)b+B}{\sqrt{(u_b/u_B)^2+1}}+\frac{(v_b/v_B)b+B}{\sqrt{(v_b/v_B)^2+1}}\right)
=\pm\frac{2B}{\sqrt{(u_b/u_B)^2+1}}\in V_M.
\end{align*}
That is, $(u,v)\in E_1=E$.

Conversely, suppose that $(u,v)\in E$, set 
\[
w=\frac u{\|u\|_1}+\frac v{\|v\|_1}\in V_m\quad\mbox{and}\quad W=\frac u{\|u\|_1}-\frac v{\|v\|_1}\in V_M,
\]
and observe that $w+W$ is in the direction of $u$ and $w-W$ is in the direction of $v$. By Lemma \ref{orthogonal}, $w$ and $W$ are orthogonal with respect to both inner products. Thus,
\[
\cos\varphi=\frac{\1{w+W}{w-W}}{\|w+W\|_1\|w-W\|_1}=\frac{\|w\|_1^2-\|W\|_1^2}{\|w\|_1^2+\|W\|_1^2}
\]
and
\[
\tan^2(\varphi/2)=\frac{1-\cos\varphi}{1+\cos\varphi}=\frac{\|W\|_1^2}{\|w\|_1^2}.
\]
A similar calculation yields the corresponding formula for $\psi$ and leads to the conclusion,
\[
\tan^2(\psi/2)=\frac{\|W\|_2^2}{\|w\|_2^2}=\frac{M^2\|W\|_1^2}{m^2\|w\|_1^2}=(M/m)^2\tan^2(\varphi/2).
\]
Taking square roots establishes equality in the right-hand inequality of (\ref{TAN}).

Applying the right-hand inequality of (\ref{TAN}) to the vectors $u$ and $-v$ replaces $\varphi$ by $\pi-\varphi$ and $\psi$ by $\pi-\psi$ to give the conclusion,
\[
\cot(\psi/2)=\tan(\pi/2-\psi/2)\le (M/m)\tan(\pi/2-\varphi/2)=(M/m)\cot(\varphi/2).
\]
This proves the left-hand inequality of (\ref{TAN}), with equality if and only if $(u,-v)\in E$. This completes the proof in the case that $V$ is a real vector space.

If $V$ is a complex space and $\1\cdot\cdot$ and $\2\cdot\cdot$ are complex inner products, the conclusion of the theorem follows by applying the result just proved to the real vector space $V_\mathbb R$ equipped with the real inner products $\re\1\cdot\cdot$ and $\re\2\cdot\cdot$.  This completes the proof.\qed\vskip 2ex

The angle between two subsets of $V$ is defined as an infimum of angles between pairs of vectors. The inequality (\ref{TAN}) remains valid when we take an infimum of all three terms so we have the following result. Note that since the cosine function is decreasing, the cosine of an infimum of angles is achieved by taking the supremum of their cosines.

\begin{corollary}\label{ST} Let $V$ be a real or complex vector space equipped with inner products $\1\cdot\cdot$ and $\2\cdot\cdot$. Make the definitions (\ref{defs}). For $S,T\subseteq V$, each containing at least one non-zero vector, let $\Phi$ and $\Psi$ be the angles between the subsets $S$ and $T$ with respect to $\1\cdot\cdot$ and $\2\cdot\cdot$, respectively. That is, $0\le\Phi\le\pi$, $0\le\Psi\le\pi$,
\begin{equation}\label{sups}
\cos\Phi=\sup_{\substack{0\ne u\in S\\0\ne v\in T}}\frac{\re\1uv}{\|u\|_1\|v\|_1},\quad\mbox{and}\quad
\cos\Psi=\sup_{\substack{0\ne u\in S\\0\ne v\in T}}\frac{\re\2uv}{\|u\|_2\|v\|_2}.
\end{equation}
Then
\[
(m/M)\tan(\Phi/2)\le\tan(\Psi/2)\le (M/m)\tan(\Phi/2).
\] 
\end{corollary}

The following theorem is our version of the generalized Wielandt inequality in inner product spaces. As pointed out earlier, the angles between the (real or complex) lines determined by $u$ and $v$ are often taken as alternative definitions of the angle between vectors themselves. We show that with this definition the results of Theorem \ref{main} still hold, but the conditions for equality become slightly more complicated. 

\begin{theorem}\label{lines} Let $V$ be a real or complex vector space equipped with inner products $\1\cdot\cdot$ and $\2\cdot\cdot$. Make the definitions (\ref{defs}). For independent vectors $u$ and $v$ in $V$ let $\Phi$ and $\Psi$ be defined by, $0\le\Phi\le\pi/2$, $0\le\Psi\le\pi/2$,
\[
\cos\Phi=\frac{|\1uv|}{\|u\|_1\|v\|_1}\quad\mbox{and}\quad\cos\Psi=\frac{|\2uv|}{\|u\|_2\|v\|_2}.
\]
Then 
\begin{equation}\label{TAN2}
(m/M)\tan(\Phi/2)\le\tan(\Psi/2)\le (M/m)\tan(\Phi/2).
\end{equation}
Let $\alpha_1$ and $\alpha_2$ be solutions to $|\1uv|=\alpha_1\1uv$ and $|\2uv|=\alpha_2\2uv$. Equality holds in the right-hand inequality of (\ref{TAN2}) if and only if $(\alpha_1u,v)\in E$ and either $\alpha_1=\alpha_2$ or $\2uv=0$.
Equality holds in the left-hand inequality of (\ref{TAN2}) if and only if $(\alpha_2u,-v)\in E$ and either $\alpha_1=\alpha_2$ or $\1uv=0$.
\end{theorem}
Proof. Apply Corollary \ref{ST} to the lines $S=\mathbb C u$ and $T=\mathbb C v$  ($S=\mathbb R u$ and $T=\mathbb R v$ in the real case) to obtain (\ref{TAN2}). By (\ref{alpha}), $\Phi$ is the angle between $\alpha_1u$ and $v$ with respect to $\1\cdot\cdot$ and $\Psi$ is the angle between $\alpha_2u$ and $v$ with respect to $\2\cdot\cdot$. To analyse the right-hand inequality of (\ref{TAN2}), let $\theta$ be the angle between $\alpha_1u$ and $v$ with respect to $\2\cdot\cdot$. The infimum definition of $\Psi$ and  Theorem \ref{main} show that
\begin{equation}\label{thetaright}
\tan(\Psi/2)\le\tan(\theta/2)\le(M/m)\tan(\Phi/2).
\end{equation}
By (\ref{alpha}), the first of these is equality if and only if either $\alpha_1=\alpha_2$ or $\2uv=0$. By Theorem \ref{main}, the second is equality if and only if $(\alpha_1u,v)\in E$. Thus equality holds in the right-hand inequality of (\ref{TAN2}) if and only if  $(\alpha_1u,v)\in E$ and either $\alpha_1=\alpha_2$ or $\2uv=0$.

To analyse the left-hand inequality of (\ref{TAN2}), let $\theta$ be the angle between $\alpha_2u$ and $v$ with respect to $\1\cdot\cdot$. The infimum definition of $\Phi$ and  Theorem \ref{main} show that
\begin{equation}\label{thetaleft}
(m/M)\tan(\Phi/2)\le(m/M)\tan(\theta/2)\le\tan(\Psi/2).
\end{equation}
By (\ref{alpha}), the first of these is equality if and only if either $\alpha_1=\alpha_2$ or $\1uv=0$. By Theorem \ref{main}, the second is equality if and only if $(\alpha_2u,-v)\in E$. Thus equality holds in the left-hand inequality of (\ref{TAN2}) if and only if $(\alpha_2u,-v)\in E$ and either $\alpha_1=\alpha_2$ or $\1uv=0$. \qed\vskip 2ex

The inequalities (\ref{TAN}) and (\ref{TAN2}) can be expressed in various equivalent forms. In terms of cosines (\ref{TAN}) becomes, with $\chi=(M^2-m^2)/(M^2+m^2)$, 
\begin{equation}\label{COS}
\frac{-\chi+\cos\varphi}{1-\chi\cos\varphi}
\le\cos\psi
\le\frac{\chi+\cos\varphi}{1+\chi\cos\varphi}.
\end{equation}
Replace $\varphi$ and $\psi$ by $\Phi$ and $\Psi$ to get the
expression for (\ref{TAN2}).
In terms of inner products instead of angles, the inequalities (\ref{TAN}) of Theorem \ref{main} and (\ref{TAN2}) of Theorem \ref{lines} become, in the case $\|u\|_1=\|v\|_1=1$,
\begin{equation}\label{ReIP}
\frac{-\chi+\re\1uv}{1-\chi\re\1uv}
\le\frac{\re\2uv}{\|u\|_2\|v\|_2}
\le\frac{\chi+\re\1uv}{1+\chi\re\1uv}.
\end{equation}
and
\begin{equation}\label{modIP}
\frac{-\chi+|\1uv|}{1-\chi|\1uv|}
\le\frac{|\2uv|}{\|u\|_2\|v\|_2}
\le\frac{\chi+|\1uv|}{1+\chi|\1uv|},
\end{equation}
respectively.

The special case $\Phi=\pi/2$ in Theorem \ref{lines} gives an inner product formulation of Wielandt's inequality that includes all cases of equality.  Note that the right-hand inequality of (\ref{modIP}) is equivalent to the left-hand inequality of (\ref{TAN2}).

\begin{corollary} Let $V$ be a real or complex vector space equipped with inner products $\1\cdot\cdot$ and $\2\cdot\cdot$. Make the definitions (\ref{defs}). Suppose the non-zero vectors $u,v\in V$ are orthogonal with respect to $\1\cdot\cdot$ and $\alpha$ satisfies $|\2uv|=\alpha\2uv$. Then,
\begin{equation}\label{orthIP}
\frac{|\2uv|}{\|u\|_2\|v\|_2}
\le\frac{M^2-m^2}{M^2+m^2}
\end{equation}
with equality if and only if $(\alpha u,-v)\in E$.
\end{corollary}

The following theorem gives upper and lower bounds on the difference between the cosines of $\varphi$ and $\psi$. It improves the estimates given in Theorems 1 and 2 of \cite{D}.

\begin{theorem}\label{Drag} Let $V$ be a real or complex vector space equipped with inner products $\1\cdot\cdot$ and $\2\cdot\cdot$. Make the definitions (\ref{defs}).   For independent vectors $u$ and $v$ in $V$,
\begin{equation}\label{regen}
-2\frac{M-m}{M+m}\le\frac{\re\2uv}{\|u\|_2\|v\|_2}-\frac{\re\1uv}{\|u\|_1\|v\|_1}\le2\frac{M-m}{M+m}
\end{equation}
and, if $\re\1uv\ge0$, then
\begin{equation}\label{pos}
\frac{\re\2uv}{\|u\|_2\|v\|_2}-\frac{\re\1uv}{\|u\|_1\|v\|_1}\le\frac{M^2-m^2}{M^2+m^2}.
\end{equation}
Also,
\begin{equation}\label{absgen}
-\frac{M^2-m^2}{M^2+m^2}\le\frac{|\2uv|}{\|u\|_2\|v\|_2}-\frac{|\1uv|}{\|u\|_1\|v\|_1}\le\frac{M^2-m^2}{M^2+m^2}.
\end{equation}
\end{theorem}
Proof.   Suppose $\varphi$ and $\psi$ are the angles between $u$ and $v$ with respect to $\1\cdot\cdot$ and $\2\cdot\cdot$. Since,
\[
\cos\psi-\cos\varphi=2/(1+\tan^2(\psi/2))-2/(1+\tan^2(\varphi/2)),
\]
Theorem \ref{main} gives
\[
\frac2{1+(M/m)^2x}-\frac2{1+x}\le\cos\psi-\cos\varphi\le\frac2{1+(m/M)^2x}-\frac2{1+x},
\]
where $x=\tan^2(\varphi/2)$.  A little calculus shows that the minimum value, over all $x\in[0,\infty]$, of the expression on the left occurs at $x=m/M$ and the maximum value, over all $x\in[0,\infty]$, of the expression on the right occurs at $x=M/m$. This gives (\ref{regen}). If $\re\1uv\ge0$ then $\varphi\le\pi/2$ and so $x =\tan^2(\varphi/2)\le 1$. The maximum value on the right now occurs at $x=1$, giving (\ref{pos}).

The same analysis, applied to the angles $\Phi$ and $\Psi$ between the lines $\mathbb C u$ and $\mathbb C v$ (or $\mathbb R u$ and $\mathbb R v$ in the real case) includes the restriction $\tan^2(\Phi/2)\le 1$ and gives the right-hand inequality in (\ref{absgen}). The left-hand inequality follows from the right-hand one by interchanging the inner products $\1\cdot\cdot$ and $\2\cdot\cdot$. Besides interchanging the angles $\varphi$ and $\psi$, this has the effect of replacing $m$ by $1/M$ and $M$ by $1/m$ to give
\[
\frac{|\1uv|}{\|u\|_1\|v\|_1}-\frac{|\2uv|}{\|u\|_2\|v\|_2}\le\frac{(1/m)^2-(1/M)^2}{(1/m)^2+(1/M)^2}=\frac{M^2-m^2}{M^2+m^2}.
\]
Multiplying through by $-1$ completes the proof.\qed\vskip 2ex

In our notation, Dragomir's results from \cite{D} are
\[
1-\frac{M^2}{m^2}\le\frac{|\2uv|}{\|u\|_2\|v\|_2}-\frac{|\1uv|}{\|u\|_1\|v\|_1}\le1-\frac{m^2}{M^2},
\]
and, if $\re\1uv\ge0$, then
\[
1-\frac{M^2}{m^2}\le\frac{\re\2uv}{\|u\|_2\|v\|_2}-\frac{\re\1uv}{\|u\|_1\|v\|_1}\le1-\frac{m^2}{M^2}.
\]
Since
\[
1-\frac{M^2}{m^2}\le-2\frac{M-m}{M+m}\le -\frac{M^2-m^2}{M^2+m^2} \quad\mbox{and}\quad\frac{M^2-m^2}{M^2+m^2}\le1-\frac{m^2}{M^2},
\]
Theorem \ref{Drag} improves on both of these statements.

The estimate (\ref{regen}), on the difference between the cosines of $\varphi$ and $\psi$ readily gives a lower bound on the product of those cosines.
\begin{corollary}\label{cosprod}  Let $V$ be a real or complex vector space equipped with inner products $\1\cdot\cdot$ and $\2\cdot\cdot$. Make the definitions (\ref{defs}). For independent vectors $u$ and $v$ in $V$,
\begin{equation}\label{coscos}
\frac{\re\1uv}{\|u\|_1\|v\|_1}\frac{\re\2uv}{\|u\|_2\|v\|_2}\ge-\left(\frac{M-m}{M+m}\right)^2.
\end{equation}
\end{corollary}
Proof. Let $\mu=(M-m)/(M+m)$,
\[
x=\frac{\re\1uv}{\|u\|_1\|v\|_1},\quad\mbox{and}\quad y=\frac{\re\2uv}{\|u\|_2\|v\|_2}.
\]
Note that $0\le\mu<1$. By the Cauchy-Schwarz inequality and (\ref{regen}), the point $(x,y)$ lies in the region defined by $-1\le x\le1$, $-1\le y\le1$, and $-2\mu\le x-y\le 2\mu$. Minimizing $xy$ over this hexagonal region easily yields $(x,y)=(-\mu,\mu)$ or $(x,y)=(\mu,-\mu)$. Thus, $xy\ge-\mu^2$ as required. \qed

\section{Formulation in terms of matrices}\label{matrices}

The angle $\theta$ between vectors $x,y\in \mathbb C^n$ is defined by $0\le\theta\le\pi$ and
\[
\cos\theta=\frac{\re y^*x}{\|x\|\|y\|} 
\]
and the angle $\Theta$ between the complex lines $\mathbb Cx$ and $\mathbb Cy$ satisfies $0\le\Theta\le\pi/2$ and
\[
\cos\Theta=\frac{|y^*x|}{\|x\|\|y\|}. 
\]
Let $A$ be an invertible $n\times n$ matrix and consider the two inner products 
\begin{equation}\label{MIPs}
\1xy=y^*x\quad\mbox{and}\quad\2xy=(Ay)^*(Ax)
\end{equation}
on $\mathbb C^n$. Then the definitions in (\ref{defs}) show that $M=\|A\|$ and $1/m=\|A^{-1}\|$ so the condition number of $A$ is $\kappa(A)=M/m$. Theorem \ref{main} becomes the following.
\begin{theorem}\label{Mmain} Let $A$ be an invertible $n\times n$ matrix. For independent $x,y\in\mathbb C^n$ let $\varphi$ be the angle between $x$ and $y$ and let $\psi$ be the angle between $Ax$ and $Ay$. Then, 
\[
\kappa(A)^{-1}\tan(\varphi/2)\le\tan(\psi/2)\le \kappa(A)\tan(\varphi/2).
\]
Let $\lambda_n$ and $\lambda_1$ denote the smallest and largest eigenvalues of $A^*A$. Then equality holds in the right-hand inequality above if and only if $x/\|x\|+y/\|y\|$ is in the $\lambda_n$-eigenspace of $A^*A$ and $x/\|x\|-y/\|y\|$ is in the $\lambda_1$-eigenspace of $A^*A$. Also, equality holds in the left-hand inequality above if and only if $x/\|x\|-y/\|y\|$ is in the $\lambda_n$-eigenspace of $A^*A$ and $x/\|x\|+y/\|y\|$ is in the $\lambda_1$-eigenspace of $A^*A$.
\end{theorem}
Theorem \ref{lines} gives a concise reformulation of the generalized Wielandt inequality. Since $\kappa(A)=\kappa(A^{-1})$, the symmetry between the angles $\Phi$ and $\Psi$ is clear.
\begin{theorem}\label{Mlines} Let $A$ be an invertible $n\times n$ matrix. For independent $x,y\in\mathbb C^n$ let $\Phi$ be the angle between the complex lines $\mathbb Cx$ and $\mathbb Cy$ and let $\Psi$ be the angle between the complex lines $\mathbb C(Ax)$ and $\mathbb C(Ay)$. 
Then
\[
\kappa(A)^{-1}\tan(\Phi/2)\le\tan(\Psi/2)\le \kappa(A)\tan(\Phi/2).
\]
\end{theorem}
It takes a bit of care to show the equivalence of this theorem with Theorem \ref{gW} because the angles $\Phi$ and $\Psi$ represent subtly different concepts in the two statements. In Theorem \ref{Mlines}, $\Phi$ and $\Psi$ represent angles between given complex lines, while in Theorem \ref{gW} they represent bounds on those angles rather than the angles themselves. Also, one must apply Theorem \ref{gW} to $A$ and to $A^{-1}$ (or else to $x,y$ and to $x, -y$) to obtain both sides of the inequality above. 

The conclusion of Theorems \ref{Mmain} and \ref{Mlines} may be rewritten as
\begin{equation}\label{MCOS}
\frac{-\chi+\cos\varphi}{1-\chi\cos\varphi}
\le\cos\psi
\le\frac{\chi+\cos\varphi}{1+\chi\cos\varphi},
\end{equation}
where $\chi=(\kappa(A)^2-1)/(\kappa(A)^2+1)$. (Of course, $\varphi$ and
$\psi$ should be replaced by $\Phi$ and $\Psi$ when rewriting Theorem
\ref{Mlines}.)

We have omitted the characterization of the cases of equality in Theorem \ref{Mlines} but they can be readily obtained from Theorem \ref{lines}. Conditions for equality in Theorem \ref{main} are simpler than those in Theorem \ref{lines} because the former deals with angles between a single pair of vectors and the latter with an infimum of angles between vectors in two one-dimensional subspaces. To recognize when equality occurs in Theorem \ref{main} one only has to consider the placement of the vectors $u$ and $v$ relative to the eigenspaces $V_m$ and $V_M$. But equality in Theorem \ref{lines} requires that this infimum of angles be achieved for $u$ and $v$ in addition to requiring their correct placement with respect to these eigenspaces. In \cite{K}, Kolotilina gave the following characterization of the cases of equality in the generalized Wielandt inequality, without explicit recognition of this two-stage requirement. We give an alternative proof using Theorem \ref{lines}. (Notice that the complex numbers $\xi$ and $\eta$ appearing in the Theorem of \cite{K} are unnecessary as they may be absorbed into the eigenvectors $x_1$ and $x_n$.)

\begin{proposition} Let $B$ be an $n\times n$ invertible Hermitian matrix, suppose $\lambda_1>\lambda_n>0$ are its largest and smallest eigenvalues, respectively, and set $\chi=(\lambda_1-\lambda_n)/(\lambda_1+\lambda_n)$. Fix independent $x,y\in\mathbb C^n$ and let $\cos\varphi=|y^*x|/(\|x\|\|y\|)$. Then
\begin{equation}\label{kolo}
|y^*Bx|=\frac{\chi+\cos\varphi}{1+\chi\cos\varphi}\sqrt{x^*Bx}\sqrt{y^*By}
\end{equation}
if and only if
\begin{equation}\label{kolocond}
\begin{aligned}
\frac{x}{\|x\|}&=\frac1{\sqrt2}(\sqrt{1+\cos\varphi}\,x_1+\sqrt{1-\cos\varphi}\,x_n),\quad\mbox{and}\\
\frac{y}{\|y\|}&=\frac\varepsilon{\sqrt2}(\sqrt{1+\cos\varphi}\,x_1-\sqrt{1-\cos\varphi}\,x_n)
\end{aligned}
\end{equation}
for some complex number $\varepsilon$ of unit modulus and some unit eigenvectors $x_1$ and $x_n$ satisfying $Bx_1=\lambda_1x_1$ and $Bx_n=\lambda_nx_n$.
\end{proposition}
Proof. With $A=B^{1/2}$ we have $B=A^*A$. Apply Theorem \ref{lines} to the inner products (\ref{MIPs}) and note that $M=\lambda_1$ and $m=\lambda_n$ so $V_M$ and $V_m$ are the $\lambda_1$- and $\lambda_n$-eigenspaces of $B$, respectively. Using (\ref{COS}), we see that (\ref{kolo}) is equivalent to equality in the left hand inequality of (\ref{TAN2}). Thus, Theorem \ref{lines} shows that (\ref{kolo}) holds if and only if $(\alpha_2x,-y)\in E$ and either $\alpha_1=\alpha_2$ or $y^*x=0$. As in Theorem \ref{lines},   $|y^*x|=\alpha_1y^*x$ and $|(Ay)^*(Ax)|=\alpha_2(Ay)^*(Ax)$.

First suppose that  $x$ and $y$ satisfy (\ref{kolocond}). A calculation, using the fact that $x_1$ and $x_n$ are orthogonal, shows that $\varepsilon y^*x\ge0$ and $\varepsilon (Ay)^*(Ax)\ge0$. It follows that either $\alpha_1=\alpha_2=\varepsilon$ or $y^*x=0$. Also,
\[
\frac{\varepsilon x}{\|\varepsilon x\|}+\frac{-y}{\|-y\|}=\sqrt2\varepsilon\sqrt{1-\cos\varphi}\,x_n\in V_m
\]
and
\[
\frac{\varepsilon x}{\|\varepsilon x\|}-\frac{-y}{\|-y\|}=\sqrt2\varepsilon\sqrt{1+\cos\varphi}\,x_1\in V_M
\]
so $(\alpha_2x,-y)\in E$.

Conversely, suppose that $(\alpha_2x,-y)\in E$ and either $\alpha_1=\alpha_2$ or $y^*x=0$. Set $\varepsilon=\alpha_2$. Then there exist $w\in V_m$ and $W\in V_M$ such that
\[
\frac{\varepsilon x}{\|x\|}-\frac{y}{\|y\|}=w\quad\mbox{and}\quad\frac{\varepsilon x}{\|x\|}+\frac{y}{\|y\|}=W.
\]
Since $w$ and $W$ are orthogonal,the parallelogram law gives $\|W\|^2+\|w\|^2=4$ and the definition of $\varphi$ gives $\|W\|^2-\|w\|^2=4\cos\varphi$. Solving these two equations yields, $\|W\|=\sqrt2\sqrt{1+\cos\varphi}$ and $\|w\|=\sqrt2\sqrt{1-\cos\varphi}$. With $x_1=\bar\varepsilon W/\|W\|$ and $x_n=\bar\varepsilon w/\|w\|$ we have (\ref{kolocond}). This completes the proof.\qed

In Theorem 3 of \cite{Yeh}, Yeh gave a different generalization of the Wielandt inequality for angles between complex lines. Here we show that Theorem \ref{Mlines} gives the stronger inequality.

\begin{theorem} \cite{Yeh} Let $A$ be an invertible $n\times n$ matrix. For independent $x,y\in\mathbb C^n$ let $\Phi$ be the angle between the complex lines $\mathbb Cx$ and $\mathbb Cy$ and let $\Psi$ be the angle between the complex lines $\mathbb C(Ax)$ and $\mathbb C(Ay)$. Define $\theta$ by  $0\le\theta\le\pi/2$ and $\cot(\theta/2)=\kappa(A)$.
If $\cos\Phi\le 1/\kappa(A)^2$, then
\begin{equation}\label{yeh}
\cos\Psi\le\cos\theta+ 2\cos^2(\theta/2) \cos\Phi.
\end{equation}
\end{theorem}
Proof. By Theorem \ref{Mlines} and (\ref{MCOS}), it is enough to show that
\[
\frac{\chi+\cos\Phi}{1+\chi\cos\Phi}\le\cos\theta+(1+\cos\theta)\cos\Phi,
\]
where
\[
\chi=\frac{\kappa(A)^2-1}{\kappa(A)^2+1}
=\frac{\cot^2(\theta/2)-1}{\cot^2(\theta/2)+1}=\cos\theta.
\]
But both $\chi$ and $\cos\Phi$ are positive, so
\[
\frac{\chi+\cos\Phi}{1+\chi\cos\Phi}\le\chi+\cos\Phi\le\chi+(1+\chi)\cos\Phi
\]
as required.\qed

In Theorem 3.1 of \cite{Yan}, Yan generalized the Wielandt inequality for real symmetric matrices as follows.
\begin{theorem}\label{yan}\cite{Yan}  Let $B$ be a real $n\times n$ symmetric positive definite matrix with  eigenvalues $\lambda_1\ge\lambda_2\dots \ge \lambda_n>0$. For independent $x,y\in\mathbb R^n$ define $\Phi$ by $0\le\Phi\le\pi/2$ and $\|x\|\|y\|\cos\Phi=|y^Tx|$. Then,
\begin{equation}\label{yanr}
|x^TBy|\le   \left(\max_{i,j}\frac{\lambda_i\cos^2(\Phi/2)-\lambda_j\sin^2(\Phi/2)}{\lambda_i\cos^2(\Phi/2)+\lambda_j\sin^2(\Phi/2)}\right)\sqrt{x^TBx}\sqrt{y^TBy}.
\end{equation}
\end{theorem}
It was left as a conjecture in \cite{Yan} that the theorem remains true for complex vectors $x$ and $y$ and a positive definite Hermitian matrix $B$. 

It is routine to verify that the expression
\[
\frac{s\cos^2(\Phi/2)-t\sin^2(\Phi/2)}{s\cos^2(\Phi/2)+t\sin^2(\Phi/2)}
\]
is increasing in $s$ and decreasing in $t$. Thus, the maximum in (\ref{yanr}) is achieved when $i=1$ and $j=n$, where it takes the value,
\[
\frac{\lambda_1\cos^2(\Phi/2)-\lambda_n\sin^2(\Phi/2)}{\lambda_1\cos^2(\Phi/2)+\lambda_n\sin^2(\Phi/2)}
=\frac{\chi+\cos\Phi}{1+\chi\cos\Phi}.
\]
Here $\chi=(\lambda_1/\lambda_n-1)/(\lambda_1/\lambda_n+1)$. If $A=B^{1/2}$, then $\kappa(A)^2=\kappa(B)=\lambda_1/\lambda_n$ so
Theorem \ref{Mlines} and (\ref{MCOS}) implies that Theorem \ref{yan} holds in both the real and complex cases, confirming Yan's conjecture.

We end with an improvement of Lemma 2.2 from \cite{LS}. It follows directly from Corollary \ref{cosprod} with $\1xy=y^TAx$ and $\2xy=y^TBx$.

\begin{lemma} Suppose $A$ and $B$ are real symmetric positive definite $n\times n$ matrices and let $\kappa=\kappa(A^{-1/2}BA^{-1/2})$. Then for $x,y\in \mathbb R^n$ with $y\ne0$,
$$
\frac{y^TAx}{\sqrt{x^TAx}\sqrt{y^TAy}}\frac{y^TBx}{\sqrt{x^TBx}\sqrt{y^TBy}}\ge-\left(\frac{\sqrt\kappa-1}{\sqrt\kappa+1}\right)^2.
$$
\end{lemma}
The above inequality followed by the AM-GM inequality give the conclusion of Lemma 2.2 from \cite{LS}:
\begin{align*}
2\frac{y^TAx}{y^TAy}\frac{y^TBx}{y^TBy}
&\ge-2\left(\frac{\sqrt\kappa-1}{\sqrt\kappa+1}\right)^2\left(\frac{x^TAx}{y^TAy}\frac{x^TBx}{y^TBy}\right)^{1/2}\\
&\ge-\left(\frac{\sqrt\kappa-1}{\sqrt\kappa+1}\right)^2\left(\frac{x^TAx}{y^TAy}+\frac{x^TBx}{y^TBy}\right).
\end{align*}

\end{document}